\documentclass[notitlepage]{article}
\usepackage{amsmath, amssymb, graphicx, geometry}
\newcommand{\hf}{\frac{1}{2}}
\renewcommand{\vec}[1]{\mathbf{#1}}
\newcommand{\mat}[1]{\mathbf{#1}}
\newcommand{\RN}{{\mathbb{R}^N}}

\begin{document}
\title{Numerical integration of discontinuous functions in many dimensions}

\author{Vassilis Pandis}

\maketitle
\begin{abstract}
We consider the problem of numerically integrating functions with hyperplane discontinuities 
over the entire Euclidean space in many dimensions. We describe a simple process through which the 
Euclidean space is partitioned into simplices on which the integrand is smooth, generalising 
the standard practice of dividing the interval used in one-dimensional problems. Our procedure is 
combined with existing adaptive cubature algorithms to significantly reduce the necessary number
of function evaluations and memory requirements of the integrator. The method is embarrassingly 
parallel and can be trivially scaled across many cores with virtually no overhead.
Our method is particularly pertinent to the integration of Green's functions, a problem directly related
to the perturbation theory of impurity models. In three spatial dimensions we observe
a speed-up of order $100$ which increases with increasing dimensionality.
\end{abstract}

\section{Introduction}\label{sec:intro}

The numerical integration of functions in many dimensions has been a central topic in numerical
analysis for a long time. Current schemes such as adaptive cubature and adaptive Monte Carlo 
perform best for smooth integrand functions. However integrands with discontinuities can arise
quite naturally in a variety of contexts, such as the calculation of the fermionic self-energy
in condensed-matter physics or the study of multiphase flows in the context of computational fluid dynamics.  

In one dimension discontinuities are easily accommodated within an adaptive framework simply
by dividing the region of integration into sub-regions on which the integrand is smooth. 
In this paper we show how this process can be extended to higher dimensionalities. Our method is applicable
to integrals which are discontinuous on any number of hyperplanes that contain the origin, and in any number of 
dimensions. We limit our attention to integrals over the entire $\RN$ --- this
is not a material limitation as integrals over a proper hyperrectangle can be straightforwardly mapped onto
$\RN$. We assume that the discontinuities in question arise from terms of the form sign$(C_{\vec{x}})$ where
$C_{\vec{x}}$ is any linear combination of the coordinates. These is precisely the form of the discontinuities 
encountered in the Green's functions of fermionic systems. 

We write the integral in question as 
\begin{equation}
I = \int_\RN   \prod_{i=1}^{M} F_i(\vec{x})d\vec{x}, 
\label{eq:integral}
\end{equation}
where $F_i(\vec{x})$ is discontinuous on the hyperplane with equation $\vec{a_i}\cdot \vec{x} = 0$. 
We construct the $M\times N$ \emph{discontinuity} matrix $\mat{C}$ such that $C_{ij} = (\vec{a_i})_j$. 
We will here assume $M \geq N$ --- we will comment on this at end the next section. 

\section{Method}

Let $S$ be the set of all $M\times M$ diagonal matrices with diagonal components $\pm 1$ ($|S| = 2^M$).  
To determine the regions on which the integrand is continuous we thus have to solve the homogeneous system of simultaneous
inequalities 
\begin{equation}
\mat{C}\mat{S}_i\vec{x} \geq 0,
\label{eq:system}
\end{equation}
for every $\mat{S}_i \in S$. Each inequality defines a closed-half space; the solution to the system is 
the intersection of these half-spaces which can be interpreted geometrically as a convex polytope in 
its half-space representation  (see \cite[p. 31]{Grunbaum67}). In the case of a homogeneous system the resultant polytope is 
in fact a polyhedral (infinite) convex cone~\cite{skeleton}. 

Let $P$ be the set of cones obtained by solving Eq.~\eqref{eq:system}.
A set of vectors $W_K = \left\{\vec{w}_1, \vec{w}_2, \ldots ..\vec{w}_p \right\}$ is a \emph{skeleton} of a cone $K$
if $\vec{x} = \sum_{i=1}^p \lambda_i \vec{w_i}$ belongs to $K$ for every $\lambda_i \geq 0$~\cite{skeleton}. 
The duality in the representation of a cone as either a system of linear inequalities or a conical combination 
of the skeleton is the essence of the well known Weyl-Minkowski theorem on cones. As we are only interested in 
subspaces of $\RN$ with dimension $N$ --- lower-dimensional subspaces correspond to polyhedral facets which do 
not contribute to the integral --- we can assume that $p \geq N$. 
The skeleton of an acute cone is unique up to scalar multiplication of the vectors~\cite{skeleton}. 
Once the normalization of the skeleton is fixed, each point in $K$ can be specified through its $\lambda$ coefficients $(\lambda_1, \ldots, \lambda_p)$.

Having achieved our goal of partitioning $\RN$ into regions where the integrand is 
continuous we now have to consider how to perform the integration over a cone $K\in P$. When $p=N$ the polytope 
constitutes an $N$-simplex which can be readily mapped onto the positive orthant by exploiting the bijection between
$\vec{x} \in K$ and $\vec{\lambda}$. When $p>N$ the situation is more complex, for the skeleton is linearly dependent
and there is no bijection to be exploited. To overcome this problem, each cone $K$ is decomposed into $N$-simplices
$\gamma^{K}_1, \gamma^{K}_2., \ldots$ which can then be individually mapped onto the positive orthant. The
set of all simplices $\digamma = \left\{ \gamma^K_1, \gamma^K_2, \ldots |  K \in P \right\}$ evidently partitions $\RN$; the original 
integration problem has thus been broken down into multiple, separate integrations, one over each 
simplex in $\digamma$. The method is inherently parallel --- barring error control considerations each region 
of integration can be processed independently of the others.

To control the precision of the calculation we use an unsophisticated two-pass scheme. The first pass consists of a 
crude integration over every simplex $\gamma \in \digamma$, with a relative precision of $10\%$, yielding a
result $\mu^{(1)}_\gamma$ with an associated error $\sigma^{(1)}_\gamma$. From the $\mu^{(1)}_\gamma$ we determine the 
simplex which contributes the most; let $\mu_{\textrm{max}} = \textrm{max} \left\{|\mu_\gamma|, \gamma \in \digamma \right\}$. 
To achieve a requested relative precision $f$ on the entire integral $I$ we then repeat the integration, now evaluating 
each simplex to an \emph{absolute} precision given by $\epsilon_{\textrm{abs}} = f\mu_\textrm{max}/\sqrt\nu$, where $\nu = |\digamma|$
denotes the total number of simplices, ensuring obviously we do not re-evaluate the regions for which $\sigma^{(1)}_\gamma < \epsilon_\textrm{abs}$.
The end result $I=\sum \mu^{(2)}_\gamma$ is then associated with an absolute error
\begin{equation}
\sigma_I = \sqrt{\sum_{\gamma \in \digamma} (\sigma^{(2)}_\gamma)^2/\nu}.
\end{equation}
In practice small deviations of the resultant precision for the requested precision may occur when there are
significant cancellations. This is not a particularly grave disadvantage as the actual error is always known.

Finally, we return to the question of the number of constraints. We have been assuming that the number of rows $M$ 
of the constraint matrix $\mat{C}$ is larger than the dimension of the integral, ignoring the case of an integrand which
has discontinuities on fewer than $N$ planes. This is dealt with by padding the rows of $\mat{C}$ with arbitrary vectors 
(so long as they are not parallel to any other vectors) until $M\geq N$. This trick has the disadvantage of causing 
unnecessary divisions of the region of integration but is necessary to guarantee the existence of cones. 

\section{Implementation}

The process outlined above is implemented in \verb|C++| with support for matrices provided by GSL~\cite{gsl}. The
input is the integrand and the matrix $\mat{C}$ construct as above specifying the discontinuities. The first step is 
the decomposition of $\RN$ into the polyhedral cones $P$. To this end we use \verb|skeleton|~\cite{skeleton}
which implements a modified version of the Motzkin-Burger algorithm. This package is called from our
code and returns the vectors comprising the skeleton of the polyhedral cones in $P$. 

To cut the polyhedral cones in $K$ into $N$-simplices we first project the vectors in $W_K$ onto the
cone's $(N-1)$-dimensional base. By `base' here we mean the subspace obtained by subtracting from
all $\vec{w} \in W_k$ their components along the axis of the cone and then expressing them as linear
combinations of $(N-1)$ orthonormal vectors. We can then construct the desired decomposition of 
$K$ into $N$-simplices by triangulating the points in the $(N-1)$-dimensional base and then adding
the origin to these $(N-1)$-simplices. In general this triangulation is not unique. There are
several algorithms to handle the triangulation of the base. We use the Quickhull algorithm
implemented in Qhull~\cite{qhull}. 

Each point $\vec{x}$ of the simplex can be written as a conical combination of the $(\vec{\lambda})$ 
and its skeleton vectors. To map the positive orthant onto the unit hypercube we use the rule
$\lambda_i = 1/u_i - 1$. Depending on the integrand other rules may be more suitable but this
was chosen for its simplicity. 

The final step is the integration itself. We use \verb|HIntLib|~\cite{hintlib1,hintlib2}, a sophisticated
\verb|C++| library that among other things implements adaptive cubature with a variety of rules and a range of Monte
Carlo methods. It would be perhaps more efficient to use an adaptive code that can directly handle the 
simplicial geometry, such as \verb|CUBPACK|~\cite{cubpack1,cubpack2} but for practical reasons this 
approach was not followed here. The integrations are performed in parallel using $\verb|OpenMP|$ (\verb|HIntLib|'s
native parallelization is not used).

\section{Results}

We test the method with a variety of integrands and for various dimensionalities. To do so we also have
to prescribe the discontinuities. To streamline the discussion we express Eq.~\eqref{eq:integral} as
\begin{equation}
I = \int_\RN   \prod_{i=1}^{M} F(g_i(\vec{x}))d\vec{x}, 
\label{eq:integral2}
\end{equation}
where $\vec{g}(\vec{x}) = \mat{C}\vec{x}$. A variety of test-matrices $\mat{C}$ are considered --- they are listed in the Appendix. 
All integrations are done using \verb|HIntLib|'s adaptive routines and its implementation of the embedded degree-$7$ rule of Genz and Malik~\cite{Genz80}. 
We define the following integrand test-functions
\begin{align}
F_1 (u) &= \frac{1}{u - \alpha + i \beta \textrm{sign}(u)}   \label{eq:f1} \\
F_2 (u) &= \frac{1}{u^2 - \alpha + i \beta \textrm{sign}(u)} \label{eq:f2}.
\end{align}
We note that $F_1$ is actually the non-interacting Green's function for the Anderson impurity model
in the flat-band approximation~\cite{Hewson97}. Our method was developed with this integrand
in mind --- we also consider the integrand $F_{2}$ to illustrate the more general applicability of the 
method. As $F_{1}, F_{2}$ are complex-valued, we consider for brevity only the real part of Eq.~\eqref{eq:integral2}. 

We compare the speed-up afforded by our partitioning in terms of the number of integrand evaluations 
required to achieve a given relative error. In doing so, we seemingly ignore the computational effort 
required for the partitioning itself. In practice this turns out to be essentially insignificant, owing
to the large computational cost of the integrations. Nevertheless the time spent in partitioning
can be reduced by noting that the solution of Eq.~\eqref{eq:system} and subsequent triangulation can be
carried out in parallel for each $\mat{S_i} \in S$.

Unless otherwise stated all integrations are carried out to $\epsilon_{\textrm{rel}} \approx 10^{-4}$. Due to the
two-pass technique for controlling the precision of the integration it may happen that the resultant error
is (sometimes significantly) less than requested. This is to be expected when the $\mu_\gamma$ have mostly the
same sign, i.e.\ $|I| \gg \mu_{\textrm{max}}$; in such cases the target absolute error --- which is based on $\mu_{\textrm{max}}$ ---
is smaller than necessary. The resultant error may be less than the requested in
another way: To obtain an estimate over each simplex $\gamma$, \verb|HIntLib| requires a minimum number $N_\textrm{min}$
of integrand evaluations. When this yields an estimate of the integral more precise than requested, nothing
can be done to reduce the number of evaluations.

\begin{table}
\caption{Results for $F_1$, $\alpha = -0.2,\beta=0.1$. \label{table:f1}}{
\begin{tabular}{l | c | c | c | r}
$N$   &   $M$    &   $N_{p}$   &   $N_{H}$   &  $N_{p}/N$ \\
\hline
2     &   3      & $5.3 \times 10^4$   &  $2.2 \times 10^6$ &  $41.4$ \\
2     &   4      & $8.6 \times 10^4$   &  $8.3 \times 10^5$ &  $9.6$  \\
2     &   5      & $1.3 \times 10^5$   &  $1.2 \times 10^6$ &  $9.0$  \\
2     &   6      & $1.7 \times 10^5$   &  $1.0 \times 10^6$ &  $6.3$  \\
3     &   5      & $3.3 \times 10^6$   &  $>3.0\times 10^9$ &  $>906.7$ \\
3     &   6      & $7.2 \times 10^6$   &  $1.6 \times 10^9$ &  $224.9$ \\
3     &   7      & $9.6 \times 10^6$   &  $2.6 \times 10^9$ &  $265.5$ \\
3     &   8      & $1.4 \times 10^7$   &  $2.5 \times 10^9$ &  $172.9$ \\
3     &   9      & $2.0 \times 10^7$   &  $2.7 \times 10^9$ &  $132.6$ \\ 
4     &   7      & $5.8 \times 10^8$   &  $>3.0\times 10^9$ &  $>5.2$ \\
5     &   9      & $4.3 \times 10^{10}$   & - &  -
\end{tabular}}
\end{table}
\begin{table}
\caption{Results for $F_2$, $\alpha = -0.2, \beta = 0.1$.\label{table:f2}}{
\begin{tabular}{l | c | c | c | c | c | c|  r}
$N$   &   $M$    &   $N_{p}$           & $\epsilon^{*}_\textrm{rel}$    &   $N_{H}$   & $N^{*}_{H}$  &  $N_{p}/N$ & $N^{*}_{p}/N$   \\
\hline
2     &    3     &   $3.6\times 10^4$   & $1.8\times10^{-7}$   & $1.3  \times 10^4$        & $2.5 \times 10^5$  &  $0.35$      & $7.0$ \\
2     &    4     &   $4.8\times 10^4$   & $4.6\times10^{-6}$   & $4.0  \times 10^5$        & $8.1 \times 10^6$  &  $8.3$       & $169$ \\
2     &    5     &   $6.0\times 10^4$   & $4.6\times10^{-6}$   & $6.5  \times 10^5$        & $1.4 \times 10^7$  &  $10.8$      & $236$ \\
2     &    6     &   $7.2\times 10^4$   & $1.6\times10^{-6}$   & $1.8  \times 10^6$        & $1.2 \times 10^8$  &  $24.9$      & $161$ \\
3     &    5     &   $1.7\times 10^5$   & $8.8\times10^{-5}$   & $9.8  \times 10^8$        & $1.2 \times 10^9$  &  $5820$      & $7418$   \\
3     &    6     &   $2.5\times 10^5$   & $7.7\times10^{-5}$   & $1.7  \times 10^8$        & $2.7 \times 10^9$  &  $6750$      & $10660$  \\
3     &    7     &   $3.8\times 10^5$   & $6.1\times10^{-5}$   & $>3.0 \times 10^9$        & $>3.0 \times 10^9$ &  $7804$      & $>7804$ \\
3     &    8     &   $4.8\times 10^5$   & $8.3\times10^{-5}$   & $>3.0 \times 10^9$        & $>3.0 \times 10^9$ &  $6244$      & $>6244$ \\
3     &    9     &   $5.8\times 10^5$   & $7.0\times10^{-5}$   & $>3.0 \times 10^9$        & $>3.0 \times 10^9$ &  $>5172$     & $>5172$ \\
4     &    7     &   $2.7\times 10^6$   & $9.8\times10^{-5}$   & $>3.0 \times 10^9$        & $>3.0 \times 10^9$ &  $>1111$     & $>5172$ \\
5     &    9     &   $4.8\times 10^{7}$& $9.9\times10^{-5}$   & $-$                       &  $-$               &   $-$          & $-$
\end{tabular}}
\end{table}

Our results are presented in Tables~\ref{table:f1}, ~\ref{table:f2}. In both tables $N$ denotes the dimension of integration, $M$ the number of hyperplane discontinuities, 
$N_p$ the number of function evaluations required using the partitioning technique and $N_H$ the number of function
evaluations required to achieve the requested precision without utilising our partitioning scheme. To make the comparison
meaningful we obtain $N_{H}$ using the same adaptive cubature routines in \verb|HIntLib| with the embedded degree-7 Genz-Malik rule that were employed to carry out the simplicial integrations, 
with each point $\vec{x}  \in \RN$ being mapped to $\vec{t} \in [-1, 1]^N$ through $x_i = t_i/(1-t^2_i)$. We emphasise
that our proposed integration method can be used in conjunction with any integration algorithm and is not tied to this
specific Genz-Malik rule.

We note that for many of the integrations in Table~\ref{table:f2} 
we were unable to reduce the precision below a certain level. Thus for each integrand we report the relative precision $\epsilon^{*}_\textrm{rel}$ \emph{actually} reached
by our partitioning scheme. We feel it is not clear whether it would be  fairer to judge the efficacy of our method method
by comparing $N_p$ to the function evaluations required to achieve the target relative precision of $10^{-4}$ or the `accidental'
precision $\epsilon^{*}_\textrm{rel}$. We thus report both quantities, the latter denoted by $N^{*}_{H}$. From our
results it's obvious that even when calculating the integral to a precision much greater than required our method 
greatly reduces the samples required of the integrand.

For the integrations attempted without partitioning $\RN$ we had to impose a maximum of $3\times10^9$ integrand evaluations to prevent the integrator
from exhausting the $8$~GB of RAM we had at our disposal. Apart from requiring fewer integrand samples, our partitioning
method also drastically reduces the amount of memory required. This is because each the grid for each simplex can be discarded
after the integration is complete rather than having to concurrently store data for all previous grid refinements. We have
however refrained from trying to quantify the improvement in the memory requirements as this is sensitive to the
details of our implementation, our choice of integration routines and the benchmarking itself rather non-trivial, given the parallel nature
of the program. Nevertheless in Table~\ref{table:f2} it is evident that the integration becomes unmanageable without
our method even in only four dimensions. 

As the number of simplices into which $\RN$ is partitioned increases very rapidly with $N$, the success of the method depends on whether the advantages
of a smooth integrand outweigh the cost of having to set up a new adaptive grid for each simplex, and the 
`unnecessary' function evaluations due to the crudeness of our error management. It is is evident that it does;
in Table~\ref{table:f1} we see that partitioning reduces the number of required integrand evaluations
by $1-3$ orders of magnitude and in Table~\ref{table:f2} by up to $4$ orders of magnitude.

\section{Conclusions}

We have described a method enabling the numerical integration of functions featuring hyperplane discontinuities
with existing adaptive cubature schemes. We showed how to construct a set $P$ of convex polyhedral cones
that partition $\RN$. Each cone $K \in P$ is then partitioned into simplices $\gamma^{K}_1,\gamma^K_2, \ldots$
which comprise the set $\digamma$ which partitions $\RN$. Each simplex $\gamma \in \digamma$ can then
be mapped onto a hypercube and integrated using any existing multidimensional numerical integration algorithm. 
We adopted a two-pass scheme to control the precision of our calculation. This allowed the essentially
complete paralellization of the integrations. 

Our method can dramatically accelerate the evaluation of such multidimensional integrals. The reduction
in the number of integrand samples required to obtain an estimate for the integral becomes more pronounced
as the dimensionality $N$ increases, and is reduced by several orders of magnitude compared to a naive integration.
Memory requirements are also greatly improved, allowing the evaluation of integrands in higher dimensions 
than would be otherwise possible. 

The method can be improved by coupling it directly to an integrator aware of the underlying simplicial 
geometry, thereby eliminating the need for a simplex-hypercube mapping. Our precision control can also be potentially replaced with a more 
advanced scheme in which no integrand evaluations are discarded and the threads dynamically synchronised. 
This could improve the performance of our method but coordinating the threads will require some effort
programming-wise.

The author would like to thank Alex C Hewson, Paul Carter and Ioannis Pesmazoglou for helpful discussions. This work was
generously supported by the Engineering and Physical Sciences Research Council.

\section*{APPENDIX}
We list here the test matrices $\mat{C}_{c \times N}$ pertinent to Eq.~\eqref{eq:integral2}.
\begin{equation*}
C_{3\times2}=\begin{pmatrix}
1 & 0  \\
0 & 1  \\
1 & 1 
\end{pmatrix} \,
C_{4\times2}=\begin{pmatrix}
1 & 0  \\
0 & 1  \\
2 & 1  \\
1 & -1 
\end{pmatrix} \,
C_{5\times2}=\begin{pmatrix}
1 & 0  \\
0 & 1  \\
2 & 1  \\
1 & -1  \\
-1 & 2 
\end{pmatrix} \,
C_{6\times2}=\begin{pmatrix}
1 & 0  \\
0 & 1  \\
2 & 1  \\
1 & 1  \\
1 & -1 \\
-1 & 2 
\end{pmatrix}
\end{equation*}

\begin{equation*}
C_{5\times3} = \begin{pmatrix}
1  &  0   &   0 \\
0  &  1   &   0 \\
0  &  0   &   1 \\
1  &  1   &   -1 \\
-1  &  2   &   1 
\end{pmatrix} \,
C_{6\times3} = \begin{pmatrix}
1  &  0   &   0 \\
0  &  1   &   0 \\
0  &  0   &   1 \\
1  &  -1  &   1 \\
1  &  1   &  -1 \\
-1  &  2   &   1 
\end{pmatrix} \,
C_{7\times3} = \begin{pmatrix}
1  &  0   &   0 \\
0  &  1   &   0 \\
0  &  0   &   1 \\
2  &  1   &  -1 \\
1  &  1   &   1 \\
1  &  1   &  -1 \\
-1 &  \hf &  2 
\end{pmatrix}
\end{equation*}
\begin{equation*}
C_{8\times3} = \begin{pmatrix}
1  &  0   &   0 \\
0  &  1   &   0 \\
0  &  0   &   1 \\
2  &  1   &  -1 \\
1  &  1   &   1 \\
1  &  1   &  -1 \\
\hf & \hf &   1 \\
-1 &  \hf &   2 
\end{pmatrix} \,
C_{9\times3} = \begin{pmatrix}
1  &  0   &   0 \\
0  &  1   &   0 \\
0  &  0   &   1 \\
2  &  1   &  -1 \\
\hf&  2   &   1 \\
-1 &  1   &  \hf \\
2  &  -1  &  1 \\
1  &   1  &  -1 \\
-1 &  \hf &   2 
\end{pmatrix}
\end{equation*}
\begin{equation*}
C_{7\times4} = \begin{pmatrix}
1  &  0  &  0  & 0 \\ 
0  &  1  &  0  & 0 \\ 
0  &  0  &  1  & 0 \\ 
0  &  0  &  0  & 1 \\ 
1  &  1  &  1  & 1 \\ 
1  &  2  &  1  & 2 \\ 
1  &  -2 &  2  & 1 
\end{pmatrix}\,
C_{9\times5} = \begin{pmatrix}
1  &  0  &  0  & 0  & 0 \\ 
0  &  1  &  0  & 0  & 0 \\ 
0  &  0  &  1  & 0  & 0 \\ 
0  &  0  &  0  & 1  & 0 \\ 
0  &  0  &  0  & 0  & 1 \\ 

1  &  1  &  1  & 1  & 1 \\ 
\hf&  1  & \hf & 1  &\hf \\ 
-1 & -1  & \hf & 1  &2 \\ 
2  &  1  &-\hf & 2  &-\hf
\end{pmatrix}
\end{equation*}
%

\bibliographystyle{plain}
\bibliography{biblio.bib}

\begin{thebibliography}{10}

\bibitem{qhull}
C.~B. Barber, D.P. Dobkin, and H.~Huhdanpaa.
\newblock The quickhull algorithm for convex hulls.
\newblock {\em ACM Trans. Math. Softw.}, 22(4):469--483, December 1996.

\bibitem{cubpack1}
R.~Cools and A.~Haegemans.
\newblock Algorithm 824: Cubpack: a package for automatic cubature; framework
  description.
\newblock {\em ACM Trans. Math. Softw.}, 29(3):287--296, September 2003.

\bibitem{gsl}
M.~Galassi, J.~Davies, J.~Theiler, B.~Gough, G.~Jungman, M.~Booth, and
  F.~Rossi.
\newblock {\em {GNU Scientific Library: Reference Manual}}.
\newblock Network Theory Ltd., February 2003.

\bibitem{cubpack2}
A.~Genz and R.~Cools.
\newblock An adaptive numerical cubature algorithm for simplices.
\newblock {\em ACM Trans. Math. Softw.}, 29(3):297--308, September 2003.

\bibitem{Genz80}
A.~Genz and A.A. Malik.
\newblock Remarks on algorithm 006: An adaptive algorithm for numerical
  integration over an n-dimensional rectangular region.
\newblock {\em Journal of Computational and Applied Mathematics}, 6(4):295 --
  302, 1980.

\bibitem{Grunbaum67}
B.~Gr{\"u}nbaum.
\newblock {\em Convex Polytopes}.
\newblock Graduate Texts in Mathematics. Springer, 1967.

\bibitem{Hewson97}
A.C. Hewson.
\newblock {\em The Kondo problem to heavy fermions}.
\newblock Cambridge studies in magnetism. Cambridge University Press, 1997.

\bibitem{hintlib1}
R.~Sch{\"u}rer.
\newblock {\em High-dimensional numerical integration on parallel computers}.
\newblock PhD thesis, Citeseer, 2001.

\bibitem{hintlib2}
R.~Sch{\"u}rer.
\newblock Hintlib manual, 2008.

\bibitem{skeleton}
N.Yu. Zolotykh.
\newblock New modification of the double description method for constructing
  the skeleton of a polyhedral cone.
\newblock {\em Computational Mathematics and Mathematical Physics},
  52(1):146--156, 2012.

\end{thebibliography}

\end{document}